\documentclass[11pt,reqno]{amsart}

\usepackage[a4paper,margin=1.15in]{geometry}
\usepackage{amsmath,amssymb,amsfonts,amsthm,mathtools}
\usepackage{mathrsfs}
\usepackage{enumitem}
\usepackage{hyperref}
\usepackage{xcolor}
\usepackage{aliascnt}

\hypersetup{
  colorlinks=true,
  linkcolor=blue!55!black,
  citecolor=blue!55!black,
  urlcolor=blue!55!black,
  pdftitle={Common Causal Ito Flows under Nondominated Martingale Laws},
  pdfauthor={Guangqian Zhao}
}

\theoremstyle{plain}
\newtheorem{theorem}{Theorem}[section]
\newaliascnt{proposition}{theorem}
\newtheorem{proposition}[proposition]{Proposition}
\aliascntresetthe{proposition}
\newaliascnt{lemma}{theorem}
\newtheorem{lemma}[lemma]{Lemma}
\aliascntresetthe{lemma}
\newaliascnt{corollary}{theorem}

\aliascntresetthe{corollary}

\theoremstyle{definition}
\newaliascnt{definition}{theorem}
\newtheorem{definition}[definition]{Definition}
\aliascntresetthe{definition}
\newaliascnt{assumption}{theorem}
\newtheorem{assumption}[assumption]{Assumption}
\aliascntresetthe{assumption}
\newaliascnt{remark}{theorem}
\newtheorem{remark}[remark]{Remark}
\aliascntresetthe{remark}
\newaliascnt{example}{theorem}
\newtheorem{example}[example]{Example}
\aliascntresetthe{example}

\usepackage[nameinlink,capitalize]{cleveref}
\crefname{theorem}{Theorem}{Theorems}
\crefname{proposition}{Proposition}{Propositions}
\crefname{lemma}{Lemma}{Lemmas}
\crefname{corollary}{Corollary}{Corollaries}
\crefname{definition}{Definition}{Definitions}
\crefname{assumption}{Assumption}{Assumptions}
\crefname{remark}{Remark}{Remarks}
\crefname{example}{Example}{Examples}
\Crefname{theorem}{Theorem}{Theorems}
\Crefname{proposition}{Proposition}{Propositions}
\Crefname{lemma}{Lemma}{Lemmas}
\Crefname{corollary}{Corollary}{Corollaries}
\Crefname{definition}{Definition}{Definitions}
\Crefname{assumption}{Assumption}{Assumptions}
\Crefname{remark}{Remark}{Remarks}
\Crefname{example}{Example}{Examples}

\newcommand{\R}{\mathbb R}

\newcommand{\cA}{\mathcal A}
\newcommand{\cP}{\mathcal P}
\newcommand{\eps}{\varepsilon}
\newcommand{\one}{\mathbf 1}
\newcommand{\dd}{\,d}
\newcommand{\Var}{\operatorname{Var}}
\newcommand{\Argmax}{\operatorname{Argmax}}
\newcommand{\norm}[1]{\left\lVert #1\right\rVert}
\newcommand{\abs}[1]{\left\lvert #1\right\rvert}

\title[Common causal It\^o flows]
{Common Causal It\^o Flows under Nondominated Martingale Laws:\\
Capacity Cores and Robust Sensitivities}
\author{Guangqian Zhao}
\address{School of Mathematical Sciences, University of Science and Technology of China, Hefei, Anhui 230026, China}
\email{zhaoguangqian@mail.ustc.edu.cn}


\begin{document}
\raggedbottom

\begin{abstract}
Let \(\Omega=C_0([0,T];\mathbb R)\), and let \(\mathfrak M_\Lambda\) be the
laws under which the coordinate process is a continuous square-integrable
martingale satisfying \(\mathrm d\langle X\rangle_t\leq
\Lambda\,\mathrm d t\).  Regularized dyadic square sums define a single
Borel causal map \(Q:\Omega\to\mathcal A_\Lambda\) which equals the
quadratic variation almost surely under every law in \(\mathfrak M_\Lambda\).
The approximants converge uniformly in \(L^r\) at rate \(2^{-n/2}\) and
uniformly on common compact sets whose complementary upper capacities tend
to zero.

For uniformly elliptic scalar coefficients, the Lamperti transform applied
to \(Q\) produces a total Borel causal \(C^1\) flow with the cocycle
property.  Every fixed section coincides almost surely with the classical
strong It\^o solution under every model.  The flow and its initial-state Jacobian
are continuous on the same capacity cores; for differentiable
finite-dimensional coefficient families, so is the parameter tangent.  A
capacity-core \(C^1\) transfer theorem yields class-uniform \(L^p\)
Fr\'echet expansions, \(W_p\)-continuity of the first-jet laws, and
attainment for continuous field payoffs of polynomial growth.

For robust terminal payoffs, the transferred first jet gives a joint
Hadamard--Danskin formula in the coefficient parameter and initial state.
An explicit constant-volatility family exhibits an active-model switch and
a nondifferentiable robust value on an explicit switching surface.
The bracket is also invariant under continuous finite-variation
translations.  The induced response kernel, multiplied by the deterministic
volatility, agrees with the Malliavin derivative under Gaussian volatility
models.
\end{abstract}

\maketitle

\section{Introduction}\label{sec:introduction}

Let \(\Omega\) be a canonical path space and let \(\cP\) be a family of
probability laws on \(\Omega\), not assumed to be dominated.  Suppose that a
stochastic equation is well posed under every \(P\in\cP\).  The modelwise
solutions \((Y^P)_{P\in\cP}\) need not determine a single measurable map
\(Y:\Omega\to E\) such that
\[
  Y=Y^P\qquad P\text{-a.s. for every }P\in\cP.
\]
Aggregation under nondominated laws is the construction of such a common
representative; see \cite{SonerTouziZhang11,Cohen12}.

The same issue appears for a sublinear expectation, including
\(G\)-expectation \cite{Peng19}.  A representation
\[
  \mathcal E(X)=\sup_{P\in\cP}E^P[X]
\]
is used here together with a common raw representative.  A single Borel
approximation scheme is run before a law is chosen and is then identified
under every model.  Regularized dyadic square sums provide the common
bracket; the flow and both tangent fields are deterministic functionals of
that bracket and the raw driver.

A second mechanism converts deterministic regularity into probabilistic
information.  Let
\[
  c_{\cP}(A)=\sup_{P\in\cP}P(A)
\]
be the upper capacity.  Suppose that a Borel map \(F:\Omega\to E\) is
continuous on compact sets \(K_R\) satisfying
\[
  c_{\cP}(K_R^c)\longrightarrow0.
\]
The restriction \(F|_{K_R}\) is an ordinary deterministic continuous map.
The capacity estimate makes the complement negligible simultaneously for
all laws.  Together they imply continuity of the image law
\[
  P\longmapsto P\circ F^{-1}
\]
along weakly convergent sequences in \(\cP\), tightness of all output laws,
and continuity of robust expectations.  If \(\cP\) is weakly compact, the
same mechanism also gives attainment.  Related capacity and
quasi-continuity frameworks for nondominated martingale laws and
\(G\)-expectations are developed in
\cite{DenisMartini06,DenisHuPeng11,HuWangZheng16}.  The field-valued form
used below applies the same cores to an entire flow and its sensitivity
field.

For first-order transfer, let
\(F:\Omega\times U\to\mathbb B\) have a common first jet which is
continuous on \(K_R\times C\), and suppose its derivative has a robust
\(L^q\) envelope.  For every \(p<q\), core uniform continuity and a single
H\"older estimate on \(K_R^c\) yield a class-uniform \(L^p\) Fr\'echet
expansion.  They also give \(W_p\)-continuity of the derivative-inclusive
field laws.  The compact Berge--Danskin principle then identifies the
directional derivative of the worst-case value through the active models.

The application therefore reduces to constructing one family of cores on
which the full output field is continuous.  Quadratic variation is the basic
difficulty because it is not continuous in the raw uniform topology on the
whole path space.  We construct a total causal bracket functional by
regularizing completed dyadic square sums into the compact convex set of
increasing \(\Lambda\)-Lipschitz paths.  Uniform martingale estimates yield one
summable approximation scheme and hence common compact continuity cores.

The bracket aggregator is used in a scalar Lamperti equation.  The equation
is deterministic once a raw path and its aggregated bracket have been
specified.  Under each \(P\in\mathfrak M_\Lambda\), the same Borel causal
map is identified with the classical It\^o solution.  Foundational pathwise
constructions of quadratic variation and stochastic integration include
\cite{Follmer81,Bichteler81,Karandikar95,ChiuCont18}.  Common stochastic
integrals under nondominated laws are developed in
\cite{SonerTouziZhang11,Nutz12}, and SDEs for typical paths in
\cite{BartlKupperNeufeld19}.  Universal measurable solution functionals for
classical SDEs and for broad classes of semimartingale-driven SDEs are
constructed in \cite{Kallenberg96,PrzybylowiczEtAl24}.  In the
\(G\)-framework, stochastic flows and rough-path lifts appear in
\cite{Gao09,GengQianYang14}; see also the causal functional calculus of
\cite{ChiuCont22}.  General rough-path continuity is treated in
\cite{Lyons98,FrizHairer20}, while scalar reductions to pathwise ordinary
equations go back to \cite{Lamperti64,Doss77,Sussmann78}.

Here completed dyadic square sums are regularized by an explicit causal map
with a model-uniform rate.  The same approximation generates compact
capacity cores on which the bracket, flow, initial Jacobian, and coefficient
tangent are jointly continuous.  This shared first-jet regularity yields
Wasserstein continuity of the field laws and first-order sensitivity of
worst-case values.

A quantitative driver-only companion construction, with explicit capacity
tails, raw-path H\"older moduli, and stability under stopping and controlled
concatenation, was developed in \cite{ZhaoDriverOnly26}.  The present paper
uses instead the total Lipschitz-envelope bracket above and pursues a different
first-order direction: common coefficient tangents, Wasserstein stability of
first-jet laws, robust Danskin sensitivities, and Gaussian identification of
the raw response kernel.  Thus the two constructions are complementary: the
companion work emphasizes quantitative geometry of the driver cores, whereas
the present work emphasizes differentiable transfer and robust sensitivity.

When the coefficient family is differentiable in a finite-dimensional
parameter, the same deterministic equation has a parameter derivative before
any probability law is selected.  This derivative remains causal, is
continuous on the common capacity cores, and becomes the classical stochastic
parameter-variation process under every \(P\in\mathfrak M_\Lambda\).  The
resulting derivative field gives Wasserstein-law stability and a joint
\((\theta,y)\)-Danskin formula for worst-case payoffs.
Classical parameter differentiability is treated in \cite{Metivier82}.  In
the \(G\)-framework, initial-state and parameter differentiability is proved
in \cite{Lin13}, while \cite{LuoWang14} gives a deterministic ODE reduction
for scalar \(G\)-SDEs.  The entire parameter-tangent field below is aggregated
by one raw causal map and transferred through the nondominated worst-case
operation.

If \(h\) is continuous and has finite variation, the dyadic bracket satisfies
\[
  Q(x+h)=Q(x).
\]
Consequently, the common flow is differentiable in every such direction and
the derivative is governed by one deterministic response kernel.  Selecting
a Gaussian law with deterministic volatility turns finite-variation
Cameron--Martin shifts into these same directions.  In the Brownian
representation, the Malliavin derivative is the response kernel multiplied
by the deterministic volatility.  The probability law changes the
Cameron--Martin geometry; the response kernel itself is already defined on
the raw path space.  Related Gaussian and rough-path differentiability
theories may be found in \cite{FrizVictoir10,Inahama14}.

\subsection{Main results and organization}

\Cref{thm:causal-bracket} constructs the total causal bracket, its uniform
\(L^r\) approximation rate, and the common compact continuity cores.
\Cref{thm:common-flow,thm:parametric-tangent} use the same raw bracket to
construct the scalar flow, its initial Jacobian, and its coefficient tangent.
The capacity-core \(C^1\) theorem, \Cref{thm:core-C1-transfer}, transfers the
resulting first jet to uniform robust Fr\'echet expansions and Wasserstein
stability.  \Cref{thm:flow-law-transfer,thm:parametric-robust-sensitivity}
give attainment and the joint active-model Danskin formula;
\Cref{ex:active-model-switch} exhibits an explicit active-law transition.
Finally,
\Cref{thm:bv-derivative,thm:malliavin-identification} construct the driver
response kernel and identify its Gaussian Malliavin interpretation.

\Cref{sec:cores} proves the capacity-core transfer principles.
\Cref{sec:bracket} treats the bounded-volatility martingale class and its
causal bracket.  \Cref{sec:flow} constructs the common scalar flow and its
coefficient-parameter tangent.
\Cref{sec:probability-output} derives law-level and robust sensitivity
consequences.  \Cref{sec:tangent} develops finite-variation and Gaussian
tangent calculus.

\section{Capacity cores and probabilistic transfer}
\label{sec:cores}

Throughout this section, \(\Omega\) and \(E\) are Polish spaces equipped
with fixed complete compatible metrics, and
\(\cP\subset\mathfrak P(\Omega)\).  The next definition separates the
deterministic regularity of a map from the probabilistic size of the region
where that regularity is available.

For \(1\le p<\infty\), \(\mathfrak P_p(E)\) denotes the laws with finite
\(p\)-th moment relative to one, hence every, base point, and \(W_p\) is the
\(p\)-Wasserstein distance induced by the given metric on \(E\).  On Banach
field spaces we use the supremum norm.

\begin{definition}[Compact capacity core]\label{def:capacity-core}
An increasing sequence of compact sets \((K_R)_{R\ge1}\) is a compact
capacity core for \(\cP\) if
\[
  \tau_R:=c_{\cP}(K_R^c)\longrightarrow0.
\]
A Borel map \(F:\Omega\to E\) is core-continuous if \(F|_{K_R}\) is
continuous for every \(R\).
\end{definition}

The union \(G=\bigcup_RK_R\) has polar complement.  Thus a core-continuous
map is continuous on each member of an increasing compact exhaustion of one
fixed quasi-sure domain.  The same exhaustion applies simultaneously to
every model.

\subsection{Transfer on capacity cores}

\begin{proposition}[Capacity-core mapping]\label{thm:core-transfer}
Let \(F_n,F:\Omega\to E\) be Borel and let \((K_R)\) be a compact capacity
core.
\begin{enumerate}[leftmargin=2.2em,label=(\roman*)]
  \item If \(F_n\to F\) uniformly on every \(K_R\), then \(F_n\to F\) in
  capacity.
  \item If \(F\) is core-continuous and \(P_n,P\in\cP\) satisfy
  \(P_n\Rightarrow P\), then
  \(P_n\circ F^{-1}\Rightarrow P\circ F^{-1}\).  The family
  \(\{P\circ F^{-1}:P\in\cP\}\) is uniformly tight.
  \item If, in addition, \(\cP\) is weakly compact, then every
  \(\varphi\in C_b(E)\) attains its worst-case expectation over \(\cP\).
\end{enumerate}
\end{proposition}

\begin{proof}
For (i), with \(a_{n,R}=\sup_{K_R}d_E(F_n,F)\),
\[
  c_{\cP}\bigl(d_E(F_n,F)>a_{n,R}\bigr)\le\tau_R.
\]
For (ii), extend \(\varphi\circ F|_{K_R}\), for
\(\varphi\in C_b(E)\), to a function \(g_R\in C_b(\Omega)\) with the same
supremum bound.  Weak convergence and
\[
  \bigl|E^{P_n}[\varphi(F)]-E^P[\varphi(F)]\bigr|
  \le \bigl|E^{P_n}[g_R]-E^P[g_R]\bigr|
      +4\norm{\varphi}_\infty\tau_R
\]
give convergence after \(n\to\infty\) and then \(R\to\infty\).
Uniform tightness follows from the compact sets \(F(K_R)\).  The resulting
continuity of \(P\mapsto E^P[\varphi(F)]\) on a weakly compact \(\cP\)
proves (iii).
\end{proof}

\subsection{Differentiable transfer}

The next theorem transfers corewise \(C^1\) regularity to a model-uniform
\(L^p\) expansion and records the capacity-tail error explicitly.

\begin{theorem}[Capacity-core \(C^1\) transfer]
\label{thm:core-C1-transfer}
Let \(\mathbb B\) be a separable Banach space, let \(U\subset\R^d\) be
open, and let \(F:\Omega\times U\to\mathbb B\) be jointly Borel.  Suppose
that \(z\mapsto F(x,z)\) is \(C^1\) for every \(x\), with jointly Borel
derivative \(D_zF\).  Let \((K_R)\) be a compact capacity core for \(\cP\)
and write \(\tau_R=c_{\cP}(K_R^c)\).  Fix a compact set \(C\subset U\) and
\(\rho>0\) such
that
\[
  C_\rho:=\{z\in\R^d:\operatorname{dist}(z,C)\le\rho\}
  \subset U.
\]
Assume that the first jet \((F,D_zF)\) is continuous on
\(K_R\times C_\rho\) for every \(R\).  Define
\[
  H_C(x):=\sup_{z\in C_\rho}\norm{D_zF(x,z)}
\]
and suppose, for some \(1\le p<q<\infty\), that
\begin{equation}\label{eq:C1-envelope}
  M_q:=\sup_{P\in\cP}
  \bigl(E^P[H_C^q]\bigr)^{1/q}<\infty.
\end{equation}
For \(0<\delta\le\rho\), put
\[
  \omega_R(\delta)
  :=\sup_{\substack{x\in K_R,\ z\in C,\ |u|\le\delta}}
  \norm{D_zF(x,z+u)-D_zF(x,z)}.
\]
Then \(\omega_R(\delta)\downarrow0\) as \(\delta\downarrow0\), and
\begin{align}
&\sup_{\substack{P\in\cP,\ z\in C\\0<|h|\le\delta}}
 \frac{\norm{F(\cdot,z+h)-F(\cdot,z)-D_zF(\cdot,z)h}_{L^p(P;\mathbb B)}}
 {|h|}
\notag\\
&\hspace{4cm}\le
  \omega_R(\delta)+2M_q\tau_R^{1/p-1/q}.
\label{eq:C1-transfer-bound}
\end{align}
Consequently,
\begin{equation}\label{eq:robust-Frechet}
  \lim_{\delta\downarrow0}
  \sup_{\substack{P\in\cP,\ z\in C\\0<|h|\le\delta}}
  \frac{\norm{F(\cdot,z+h)-F(\cdot,z)-D_zF(\cdot,z)h}_{L^p(P;\mathbb B)}}
  {|h|}=0.
\end{equation}

Moreover, the derivative field
\[
  \mathscr D_C(x):=\bigl(D_zF(x,z)\bigr)_{z\in C}
  \in C\bigl(C;\mathcal L(\R^d,\mathbb B)\bigr)
\]
is Borel and core-continuous.  If \(P_n,P\in\cP\) and \(P_n\Rightarrow P\),
then
\begin{equation}\label{eq:derivative-field-Wp}
  W_p\bigl(P_n\circ\mathscr D_C^{-1},
            P\circ\mathscr D_C^{-1}\bigr)\longrightarrow0.
\end{equation}
If, in addition,
\begin{equation}\label{eq:first-jet-envelope}
  \sup_{P\in\cP}E^P\!\left[
    \sup_{z\in C}
    \bigl(\norm{F(\cdot,z)}+\norm{D_zF(\cdot,z)}\bigr)^q
  \right]<\infty,
\end{equation}
then the same conclusion holds for the first-jet field
\[
  \mathscr J_C(x):=
  \bigl(F(x,z),D_zF(x,z)\bigr)_{z\in C}.
\]
In \eqref{eq:derivative-field-Wp}, and for the first jet, \(W_p\) is
computed from the supremum norm of the relevant field space.
\end{theorem}

\begin{proof}
Joint continuity on the compact set \(K_R\times C_\rho\) implies
\(\omega_R(\delta)\to0\).  The fundamental theorem of calculus in
\(\mathbb B\) gives, for \(z\in C\) and \(|h|\le\rho\),
\[
  F(x,z+h)-F(x,z)-D_zF(x,z)h
  =\int_0^1
    \bigl(D_zF(x,z+th)-D_zF(x,z)\bigr)h\,\dd t.
\]
On \(K_R\), the norm of the right-hand side is bounded by
\(\omega_R(|h|)|h|\); on \(K_R^c\), it is bounded by
\(2H_C|h|\).  H\"older's inequality yields
\[
  \norm{H_C\one_{K_R^c}}_{L^p(P)}
  \le \norm{H_C}_{L^q(P)}P(K_R^c)^{1/p-1/q}.
\]
This proves \eqref{eq:C1-transfer-bound}.  First let
\(\delta\downarrow0\), and then \(R\to\infty\), to obtain
\eqref{eq:robust-Frechet}.

The Carath\'eodory measurability theorem makes \(\mathscr D_C\) a Borel
map into the separable field space.  Uniform continuity of \(D_zF\) on
\(K_R\times C\) shows that this map is continuous on each \(K_R\).
Thus \Cref{thm:core-transfer}(ii) gives weak convergence of its laws.
The random variables \(\norm{\mathscr D_C}^p\) are uniformly integrable by
\eqref{eq:C1-envelope}, so weak convergence upgrades to
\eqref{eq:derivative-field-Wp}.  The argument for \(\mathscr J_C\) is the
same, using \eqref{eq:first-jet-envelope}.
\end{proof}

\subsection{Constructing cores from a common approximation}

\begin{theorem}[Summable common approximation]\label{thm:summable-core}
Let \((E,d_E)\) be a complete separable metric space, let \(\cP\) be
uniformly tight on \(\Omega\), and let
\(F_n:\Omega\to E\) be continuous.  Suppose that there are positive numbers
\(a_n\) such that
\[
  \sum_{n\ge1}a_n<\infty,
  \qquad
  \sum_{n\ge1}
  c_{\cP}\bigl(d_E(F_{n+1},F_n)>a_n\bigr)<\infty.
\]
Then there is a Borel map \(F:\Omega\to E\) and an increasing compact
capacity core \((K_R)\) such that \(F_n\to F\) uniformly on every \(K_R\).
In particular, \(F\) is core-continuous and \(F_n\to F\) quasi surely and in
capacity.
\end{theorem}

\begin{proof}
Let
\[
  D=\{x\in\Omega:(F_n(x))_{n\ge1}\text{ is \(d_E\)-Cauchy}\}.
\]
This is a Borel set.  Capacity Borel--Cantelli and
\(\sum_n a_n<\infty\) give \(c_{\cP}(D^c)=0\).  Completeness of \(d_E\)
defines
\[
  F(x)=\lim_{n\to\infty}F_n(x),\qquad x\in D.
\]
Extend \(F\) by a fixed point of \(E\) on \(D^c\).  The extension is Borel
because it is a pointwise limit of Borel maps on the Borel set \(D\).
Given
\(\eta_R\downarrow0\), choose a compact \(H_R\subset\Omega\) and \(N_R\)
so that the capacity of \(H_R^c\) and the tail sum of the exceptional-event
capacities are each at most \(\eta_R/2\).  Then
\[
  K_R^0
  =H_R\cap
  \bigcap_{n\ge N_R}
  \{d_E(F_{n+1},F_n)\le a_n\}
\]
is compact, has complementary capacity at most \(\eta_R\), and carries
uniform convergence by the Weierstrass test.  Replacing \(K_R^0\) by
\(K_R=\bigcup_{j\le R}K_j^0\) makes the cores increasing without changing
these conclusions.  Their capacity bound and uniform convergence imply
convergence in capacity.
\end{proof}

\begin{remark}
The theorem requires one deterministic approximation sequence shared by the
whole model class.  Model-dependent subsequences need not generate a common
quasi-sure domain or common continuity cores.
\end{remark}

\section{A causal bracket for bounded-volatility martingale laws}
\label{sec:bracket}

Let
\[
  \Omega=C_0([0,T];\R)
\]
with the uniform topology, and let \(X_t(x)=x_t\) be the coordinate process.
For \(\Lambda>0\), denote by \(\mathfrak M_\Lambda\) the set of laws
\(P\in\mathfrak P(\Omega)\) under which \(X\) is a square-integrable
martingale with respect to the usual augmentation of its natural filtration.
We take the continuous predictable version of its bracket and require, on
one \(P\)-full set,
\begin{equation}\label{eq:bracket-upper-bound}
  \langle X\rangle_t^P-\langle X\rangle_s^P
  \le\Lambda(t-s),
  \qquad 0\le s\le t\le T.
\end{equation}
Thus \(\mathfrak M_\Lambda\) is the class of scalar continuous-martingale
laws satisfying the displayed upper volatility bound.

\subsection{Model class and bracket construction}

\begin{lemma}[Weak compactness]\label{lem:model-compactness}
The set \(\mathfrak M_\Lambda\) is weakly compact in
\(\mathfrak P(\Omega)\).  Moreover, for every \(r\ge2\),
\begin{equation}\label{eq:martingale-moment}
  \sup_{P\in\mathfrak M_\Lambda}
  E^P\bigl[\abs{X_t-X_s}^r\bigr]
  \le C_r\Lambda^{r/2}\abs{t-s}^{r/2}.
\end{equation}
\end{lemma}

\begin{proof}
The moment estimate follows from the Burkholder--Davis--Gundy inequality and
\eqref{eq:bracket-upper-bound}.  The case \(r=4\), together with the
Kolmogorov tightness criterion, gives tightness on \(C_0([0,T])\).

It remains to prove closedness.  Let \(P_n\Rightarrow P\) with
\(P_n\in\mathfrak M_\Lambda\).  Uniform fourth moments imply uniform
integrability of all quadratic expressions below.  For
\(0\le s\le t\le T\) and every bounded nonnegative continuous cylinder
functional \(\phi\) depending only on the path up to time \(s\), weak
convergence gives
\[
  E^P[(X_t-X_s)\phi]=0
\]
and
\[
  E^P[(X_t^2-X_s^2-\Lambda(t-s))\phi]\le0.
\]
A monotone-class argument shows that \(X\) is a square-integrable
\(P\)-martingale and \(X_t^2-\Lambda t\) is a supermartingale.  In the
Doob--Meyer decomposition,
\[
  X_t^2-\Lambda t
  =\bigl(X_t^2-\langle X\rangle_t^P\bigr)
   +\bigl(\langle X\rangle_t^P-\Lambda t\bigr).
\]
Uniqueness of the predictable finite-variation part implies that
\(\langle X\rangle^P-\Lambda t\) is decreasing.  Thus
\eqref{eq:bracket-upper-bound} holds under \(P\), proving closedness.
\end{proof}

Let \(\pi_n=\{t_k^n=kT2^{-n}:0\le k\le2^n\}\).  Define the completed
dyadic square sum
\begin{equation}\label{eq:completed-qv-sum}
  V_t^n(x)
  =
  \sum_{k:\,t_{k+1}^n\le t}
  \bigl(x_{t_{k+1}^n}-x_{t_k^n}\bigr)^2.
\end{equation}
For fixed \(n\), the map \(x\mapsto V^n(x)\) is continuous from \(\Omega\)
to the bounded functions on \([0,T]\) equipped with the supremum norm, and
it is causal.

Let \(\cA_\Lambda\) be the compact convex set
\[
  \cA_\Lambda
  =
  \{q\in C([0,T]):q_0=0,\ 0\le q_t-q_s\le\Lambda(t-s)\}.
\]
For a nonnegative increasing function \(a\) with \(a_0=0\), set
\begin{equation}\label{eq:lipschitz-envelope}
  (\mathcal R_\Lambda a)_t
  =
  \inf_{0\le u\le t}
  \{a_u+\Lambda(t-u)\}.
\end{equation}

\begin{lemma}[Causal Lipschitz envelope]\label{lem:lipschitz-envelope}
The map \(\mathcal R_\Lambda\) takes nonnegative increasing functions into
\(\cA_\Lambda\), is causal, and satisfies
\[
  \norm{\mathcal R_\Lambda a-\mathcal R_\Lambda b}_\infty
  \le\norm{a-b}_\infty.
\]
If \(a\in\cA_\Lambda\), then \(\mathcal R_\Lambda a=a\).
\end{lemma}

\begin{proof}
Write
\[
  (\mathcal R_\Lambda a)_t
  =\Lambda t+\inf_{0\le u\le t}(a_u-\Lambda u).
\]
As \(t\) grows, the running infimum can decrease by at most
\(\Lambda\) times the elapsed time because \(a\) is increasing.  Hence the
displayed function is increasing and \(\Lambda\)-Lipschitz.  Causality is
immediate.  The inequality follows from the elementary bound for the
difference of two infima.  If \(a\in\cA_\Lambda\), then
\(a_u+\Lambda(t-u)\ge a_t\), while \(u=t\) gives equality.
\end{proof}

Put
\begin{equation}\label{eq:regularized-qv}
  q^n(x)
  =
  \mathcal R_\Lambda V^n(x)
  \in\cA_\Lambda.
\end{equation}
Each map \(q^n:\Omega\to\cA_\Lambda\) is continuous and causal.  Finally define,
for every raw path \(x\in\Omega\),
\begin{equation}\label{eq:total-bracket}
  Q_t(x)=\limsup_{n\to\infty}q_t^n(x),
  \qquad 0\le t\le T.
\end{equation}
Since all \(q^n\) are increasing and \(\Lambda\)-Lipschitz, the same is true
of their pointwise limsup.  Thus \(Q(x)\in\cA_\Lambda\) for every \(x\).
Evaluation at rational times shows that
\[
  Q:\Omega\longrightarrow\cA_\Lambda
\]
is Borel, and \eqref{eq:total-bracket} is causal at every time.

\subsection{Uniform convergence under all models}

For comparison with \eqref{eq:completed-qv-sum}, define the continuous
within-cell sum
\[
  \widetilde V_t^n(X)
  =
  \sum_{k:\,t_k^n<t}
  \bigl(
    X_{t_{k+1}^n\wedge t}-X_{t_k^n}
  \bigr)^2.
\]
If \(\underline r_n\) is the left endpoint of the dyadic cell containing
\(r\), It\^o's formula gives
\begin{equation}\label{eq:qv-martingale-error}
  \widetilde V_t^n(X)-\langle X\rangle_t^P
  =
  2\int_0^t(X_r-X_{\underline r_n})\dd X_r.
\end{equation}

\begin{lemma}[Uniform dyadic estimate]\label{lem:dyadic-estimate}
For every \(r\ge2\), there is \(C_{r,\Lambda,T}<\infty\) such that
\begin{equation}\label{eq:qv-Lr-rate}
  \sup_{P\in\mathfrak M_\Lambda}
  \left(
    E^P\left[
      \norm{q^n(X)-\langle X\rangle^P}_\infty^r
    \right]
  \right)^{1/r}
  \le C_{r,\Lambda,T}2^{-n/2}.
\end{equation}
\end{lemma}

\begin{proof}
Let \(h_n=T2^{-n}\).  By BDG, \eqref{eq:qv-martingale-error}, and
\(\dd\langle X\rangle_r^P\le\Lambda\dd r\),
\[
  \left\|
    \norm{\widetilde V^n-\langle X\rangle^P}_\infty
  \right\|_{L^r(P)}
  \le
  C_r
  \left(
    E^P\left[
      \left(
        \int_0^T
        \abs{X_r-X_{\underline r_n}}^2
        \dd\langle X\rangle_r^P
      \right)^{r/2}
    \right]
  \right)^{1/r}.
\]
Jensen's inequality and Fubini's theorem give
\begin{align*}
E^P\left[
  \left(\int_0^T
    \abs{X_r-X_{\underline r_n}}^2\dd\langle X\rangle_r^P
  \right)^{r/2}
\right]
&\le
\Lambda^{r/2}T^{r/2-1}
\int_0^T E^P\!\left[
  \abs{X_r-X_{\underline r_n}}^r
\right]\dd r\\
&\le C_{r,\Lambda,T}h_n^{r/2},
\end{align*}
where \eqref{eq:martingale-moment} is used in the last step.
Consequently,
\[
  \left\|
    \norm{\widetilde V^n-\langle X\rangle^P}_\infty
  \right\|_{L^r(P)}
  \le C_{r,\Lambda,T}h_n^{1/2}.
\]

The difference between \(\widetilde V^n\) and the completed sum \(V^n\) is
bounded by the largest squared oscillation on one dyadic cell.  Hence
\begin{align*}
  E^P\left[
    \norm{\widetilde V^n-V^n}_\infty^r
  \right]
  &\le
  \sum_{k<2^n}
  E^P\left[
    \sup_{t_k^n\le t\le t_{k+1}^n}
    \abs{X_t-X_{t_k^n}}^{2r}
  \right]\\
  &\le C_{r,\Lambda,T}2^nh_n^r
   \le C_{r,\Lambda,T}h_n^{r-1},
\end{align*}
where the maximal BDG inequality is applied on each dyadic cell.  For
\(r\ge2\), the resulting \(L^r\)-norm is bounded by a constant times
\(h_n^{1/2}\).  Finally, \Cref{lem:lipschitz-envelope} is nonexpansive and
fixes \(\langle X\rangle^P\).
This proves \eqref{eq:qv-Lr-rate}.
\end{proof}

\begin{theorem}[Total causal bracket and common cores]
\label{thm:causal-bracket}
The functional \(Q\) in \eqref{eq:total-bracket} has the following
properties.
\begin{enumerate}[leftmargin=2.2em,label=(\roman*)]
  \item It is a total Borel causal map from \(\Omega\) to
  \(\cA_\Lambda\).

  \item For every \(P\in\mathfrak M_\Lambda\),
  \begin{equation}\label{eq:bracket-identification}
    Q(X)=\langle X\rangle^P,
    \qquad P\text{-a.s.}
  \end{equation}
  Moreover, for every \(r\ge2\),
  \[
    \sup_{P\in\mathfrak M_\Lambda}
    \left(
      E^P\left[\norm{q^n-Q}_\infty^r\right]
    \right)^{1/r}
    \le C_{r,\Lambda,T}2^{-n/2}.
  \]

  \item There is an increasing compact capacity core \((K_R)\) for
  \(\mathfrak M_\Lambda\) such that \(q^n\to Q\) uniformly on every
  \(K_R\).  In particular, \(Q|_{K_R}\) is continuous in the raw uniform
  topology.
\end{enumerate}
\end{theorem}

\begin{proof}
Only (ii)--(iii) remain.  Fix \(0<\gamma<1/2\) and choose \(r\ge2\).
By \Cref{lem:dyadic-estimate} and Markov's inequality,
\[
  c_\Lambda\left(
    \norm{q^{n+1}-q^n}_\infty>2^{-\gamma n}
  \right)
  \le C_{r,\Lambda,T}2^{-r(1/2-\gamma)n},
\]
where \(c_\Lambda=c_{\mathfrak M_\Lambda}\).  The right-hand side is
summable.  Thus \((q^n)\) converges uniformly outside one polar set.  Its
limit there equals the pointwise limsup \(Q\).  Under each fixed \(P\),
\Cref{lem:dyadic-estimate} also gives
\(q^n\to\langle X\rangle^P\) in probability.  Uniqueness of limits in
probability yields \eqref{eq:bracket-identification}.  The robust
\(L^r\)-estimate follows by substituting this identification in
\eqref{eq:qv-Lr-rate}.

By \Cref{lem:model-compactness}, the family \(\mathfrak M_\Lambda\) is
uniformly tight.  The maps \(q^n:\Omega\to\cA_\Lambda\) are continuous, the
thresholds \(2^{-\gamma n}\) are summable, and the preceding capacity
probabilities are summable.  Apply \Cref{thm:summable-core}.
\end{proof}

\section{Common parameterized scalar flows}
\label{sec:flow}

Let \(\Theta\subset\R^d\) be compact.  For each \(\theta\in\Theta\),
let \(b_\theta,\sigma_\theta:\R\to\R\).  We impose the following uniform
assumption.

\begin{assumption}[Coefficient class]\label{ass:coefficients}
The map
\[
  \theta\longmapsto(b_\theta,\sigma_\theta)
\]
is continuous into \(C_b^2(\R)\times C_b^3(\R)\), and there are
\(L<\infty\) and \(0<\underline\sigma\le\overline\sigma<\infty\) such that
\[
  \sup_{\theta\in\Theta}
  \bigl(
    \norm{b_\theta}_{C_b^2}
    +\norm{\sigma_\theta}_{C_b^3}
  \bigr)
  \le L,
  \qquad
  \underline\sigma\le\sigma_\theta(y)\le\overline\sigma.
\]
\end{assumption}

The uniform \(C_b^3\)-bound supplies the spatial derivatives used in the
tangent estimates of \Cref{sec:tangent}.

We use the classical Lamperti transformation \cite{Lamperti64}.  Define the
diffeomorphism and its inverse by
\begin{equation}\label{eq:Lamperti-transform}
  F_\theta(y)=\int_0^y\frac{\dd r}{\sigma_\theta(r)},
  \qquad
  G_\theta=F_\theta^{-1},
\end{equation}
and set
\begin{equation}\label{eq:transformed-coefficients}
  \widehat b_\theta(z)
  =\frac{b_\theta(G_\theta(z))}{\sigma_\theta(G_\theta(z))},
  \qquad
  \widehat c_\theta(z)
  =-\frac12\sigma_\theta'(G_\theta(z)).
\end{equation}
These functions and their first two derivatives are uniformly bounded in
the orders needed below.

For \(x\in\Omega\), \(q\in\cA_\Lambda\), \(0\le s\le T\), and
\(y\in\R\), let \(Z^{\theta,s,y}(x,q)\) be the solution of
\begin{align}
  Z_t^{\theta,s,y}(x,q)
  ={}&F_\theta(y)+x_t-x_s
  +\int_s^t\widehat b_\theta(Z_r^{\theta,s,y})\dd r
  \notag\\
  &+\int_s^t\widehat c_\theta(Z_r^{\theta,s,y})\dd q_r,
  \qquad s\le t\le T,
  \label{eq:deterministic-Lamperti}
\end{align}
and put
\begin{equation}\label{eq:deterministic-flow-map}
  \Phi_{s,t}^\theta(x,q;y)
  =G_\theta(Z_t^{\theta,s,y}(x,q)).
\end{equation}

\begin{lemma}[Deterministic Lamperti flow]\label{lem:deterministic-flow}
Equation \eqref{eq:deterministic-Lamperti} has a unique global continuous
solution.  The map
\[
  (x,q,\theta,s,t,y)\longmapsto\Phi_{s,t}^\theta(x,q;y)
\]
is jointly continuous on
\(\Omega\times\cA_\Lambda\times\Theta\times\{s\le t\}\times\R\), locally
uniformly in \(y\).  It is causal and satisfies
\begin{equation}\label{eq:deterministic-flow-cocycle}
  \Phi_{s,t}^\theta(x,q;y)
  =
  \Phi_{u,t}^\theta
  \bigl(x,q;\Phi_{s,u}^\theta(x,q;y)\bigr),
  \qquad s\le u\le t.
\end{equation}
For fixed \((x,q,\theta,s,t)\), the map in \(y\) is an
orientation-preserving \(C^1\)-diffeomorphism, with
\begin{align}
  \partial_y\Phi_{s,t}^\theta(x,q;y)
  ={}&
  \frac{\sigma_\theta(\Phi_{s,t}^\theta(x,q;y))}
       {\sigma_\theta(y)}
  \notag\\
  &\times
  \exp\left\{
    \int_s^t\widehat b_\theta'(Z_r)\dd r
    +\int_s^t\widehat c_\theta'(Z_r)\dd q_r
  \right\}.
  \label{eq:initial-Jacobian}
\end{align}
There are constants \(0<c_*\le C_*<\infty\), depending only on
\(L,\underline\sigma,\overline\sigma,\Lambda,T\), such that
\begin{equation}\label{eq:uniform-flow-Jacobian}
  c_*\le\partial_y\Phi_{s,t}^\theta(x,q;y)\le C_*
\end{equation}
for all admissible arguments.
\end{lemma}

\begin{proof}
Existence and uniqueness follow by Picard iteration and Gronwall's inequality
for the finite measure \(\dd r+\dd q_r\).  The total mass is at most
\((1+\Lambda)T\), uniformly over \(q\in\cA_\Lambda\).  Uniform convergence
of increasing paths in \(\cA_\Lambda\) implies weak convergence of their
Stieltjes measures.

To record the stability argument, let
\((x_n,q_n,\theta_n,s_n,t_n,y_n)\to(x,q,\theta,s,t,y)\).  Extend each solution
constantly to the left of its starting time.  Gronwall's inequality gives a
uniform bound on the solutions, while convergence of \((x_n)\), the common
Lipschitz bound on \((q_n)\), and the coefficient bounds give
equicontinuity.  Every subsequence therefore has a uniformly convergent
subsequence.  Along such a subsequence,
\(\widehat c_{\theta_n}(Z^n)\) converges uniformly and
\(\dd q_n\Rightarrow\dd q\).  Since every member of \(\cA_\Lambda\) is
continuous, the limiting Stieltjes measure has no atoms at \(s\) or \(t\).
The Helly--Bray theorem, first on a fixed interval and then with
\(s_n\to s\) and \(t_n\to t\), therefore gives
\[
  \int_{s_n}^{t_n}\widehat c_{\theta_n}(Z_r^n)\dd q_n(r)
  \longrightarrow
  \int_s^t\widehat c_\theta(Z_r)\dd q(r).
\]
Indeed, the error from the uniformly convergent integrands is bounded by
\[
  \Lambda T
  \norm{\widehat c_{\theta_n}(Z^n)-\widehat c_\theta(Z)}_\infty.
\]
The remaining fixed-integrand convergence follows by approximating the
interval indicator from above and below; the endpoint errors vanish because
the Stieltjes measure of the continuous path \(q\) has no atoms.
The ordinary integral and the moving endpoint terms pass to the limit
directly.  Thus every subsequential limit solves
\eqref{eq:deterministic-Lamperti}; uniqueness identifies it with \(Z\).
The same estimates are uniform for \(y\) in compact sets, proving the stated
joint continuity.

Causality and the cocycle follow from uniqueness after restricting or
splitting the equation.  Differentiating in \(y\) in Lamperti coordinates
gives a scalar linear equation.  Solving it explicitly and using
\(G_\theta'(z)=\sigma_\theta(G_\theta(z))\) gives
\eqref{eq:initial-Jacobian}.  Uniform boundedness of the exponent and the
two-sided bound on \(\sigma_\theta\) give \eqref{eq:uniform-flow-Jacobian}.
In particular,
\[
  \Phi_{s,t}^\theta(x,q;y)-\Phi_{s,t}^\theta(x,q;0)
  =\int_0^y\partial_z\Phi_{s,t}^\theta(x,q;z)\dd z.
\]
The lower bound \(\partial_z\Phi_{s,t}^\theta\ge c_*>0\) implies that the
two sides tend to \(\pm\infty\) as \(y\to\pm\infty\).  Hence the increasing
map is onto, and the inverse function theorem makes it a global
\(C^1\)-diffeomorphism.  Compare the classical stochastic-flow theory in
\cite{Kunita90}.
\end{proof}

We now insert the total causal bracket of \Cref{thm:causal-bracket}.

\begin{definition}[Common raw flow]\label{def:common-flow}
For every \(x\in\Omega\), define
\begin{equation}\label{eq:common-flow}
  S_{s,t}^\theta(x;y)
  =\Phi_{s,t}^\theta(x,Q(x);y).
\end{equation}
\end{definition}

Since \(Q(x)\in\cA_\Lambda\) for every \(x\in\Omega\),
\eqref{eq:common-flow} is defined on the whole raw path space.

\begin{theorem}[Common Borel causal It\^o flow]\label{thm:common-flow}
Under \Cref{ass:coefficients}, the field
\[
  (x,\theta,s,t,y)\longmapsto S_{s,t}^\theta(x;y)
\]
has the following properties.
\begin{enumerate}[leftmargin=2.2em,label=(\roman*)]
  \item It is Borel and causal on the full raw path space.  It satisfies the
  cocycle identity and is an orientation-preserving \(C^1\) flow in the
  initial value.  Its initial Jacobian satisfies
  \eqref{eq:uniform-flow-Jacobian}.

  \item Under every \(P\in\mathfrak M_\Lambda\), for each fixed
  \((\theta,s,y)\), the process \(S_{s,\cdot}^\theta(X;y)\) is adapted and
  is the unique strong solution of
  \begin{equation}\label{eq:modelwise-SDE}
    Y_t
    =y+\int_s^t b_\theta(Y_r)\dd r
       +\int_s^t\sigma_\theta(Y_r)\dd X_r,
    \qquad s\le t\le T.
  \end{equation}

  \item Let \((K_R)\) be the compact capacity core from
  \Cref{thm:causal-bracket}.  For every \(m<\infty\), the field-valued map
  \begin{equation}\label{eq:flow-Jacobian-field}
    \mathscr F_m(x)
    =
    \left(
      S_{s,t}^\theta(x;y),
      \partial_yS_{s,t}^\theta(x;y)
    \right)_{\theta\in\Theta,\,s\le t,\,|y|\le m}
  \end{equation}
  is Borel from \(\Omega\) into
  \begin{equation}\label{eq:field-space}
    E_m
    =C\bigl(\Theta\times\{(s,t):s\le t\}\times[-m,m];\R^2\bigr)
  \end{equation}
  and is continuous on every \(K_R\).
\end{enumerate}
\end{theorem}

\begin{proof}
Part (i) follows by composing the Borel causal map \(x\mapsto(x,Q(x))\) with
\Cref{lem:deterministic-flow}.  At a fixed time \(t\), the resulting
value depends only on the stopped path, so it is measurable with respect to
the raw history at time \(t\).  The paths are continuous; consequently the
same representative is progressively measurable under every completed
model filtration.

Fix \(P\in\mathfrak M_\Lambda\).  By
\eqref{eq:bracket-identification},
\[
  Q(X)=\langle X\rangle^P,
  \qquad P\text{-a.s.}
\]
Apply It\^o's formula to
\[
  Y_t=G_\theta\bigl(Z_t^{\theta,s,y}(X,Q(X))\bigr).
\]
Since
\[
  G_\theta'=\sigma_\theta\circ G_\theta,
  \qquad
  G_\theta''
  =(\sigma_\theta'\sigma_\theta)\circ G_\theta,
\]
and \(\widehat c_\theta=-\sigma_\theta'\circ G_\theta/2\), the two
\(\dd\langle X\rangle^P\)-terms cancel, while
\(G_\theta'\widehat b_\theta=b_\theta\circ G_\theta\).  This proves
\eqref{eq:modelwise-SDE}.  Strong uniqueness follows from the global
Lipschitz assumptions.

For (iii), \(x\mapsto(x,Q(x))\) is continuous on \(K_R\).  Joint continuity
of the deterministic flow and of \eqref{eq:initial-Jacobian}, followed by
compactness of the parameter set in \eqref{eq:field-space}, proves
continuity of \(\mathscr F_m|_{K_R}\).  Borel measurability into the
separable function space follows from joint Borel measurability and
continuity in the field parameters.
\end{proof}

\subsection{Coefficient-parameter tangents}

The compact parameter space in \Cref{ass:coefficients} gives uniformity of
the flow field.  A differentiable coefficient parameterization produces a
second tangent structure, distinct from the path translations considered in
\Cref{sec:tangent}.

\begin{assumption}[Smooth parameterization]\label{ass:parametric}
There are an open set \(U\subset\R^d\) and a compact set
\(\Theta_1\subset U\) such that
\(\Theta\subset\operatorname{int}\Theta_1\).  The coefficient map extends
to a \(C^1\) map
\[
  U\longrightarrow C_b^2(\R)\times C_b^3(\R),
  \qquad
  \theta\longmapsto(b_\theta,\sigma_\theta),
\]
and, on \(\Theta_1\), the bounds and ellipticity in
\Cref{ass:coefficients} hold together with
\[
  \sup_{\theta\in\Theta_1}
  \left(
    \norm{D_\theta b_\theta}_{\mathcal L(\R^d,C_b^2)}
    +\norm{D_\theta\sigma_\theta}_{\mathcal L(\R^d,C_b^3)}
  \right)<\infty.
\]
\end{assumption}

For \(\vartheta\in\Theta_1\), the notation
\(F_\vartheta,G_\vartheta,\Phi^\vartheta\), and \(S^\vartheta\) refers to
the same constructions as above.  Compactness and the Banach-valued
\(C^1\) assumption make all coefficient and ellipticity constants uniform
on \(\Theta_1\); in particular, \Cref{thm:common-flow} applies with
\(\Theta_1\) in place of \(\Theta\).

For \(v\in\R^d\), place a dot over a coefficient to denote its derivative
in direction \(v\); for example,
\(\dot\sigma_\theta^v=D_\theta\sigma_\theta[v]\).  Differentiating the
Lamperti transform gives
\begin{align}
  \dot F_\theta^v(y)
  &=-\int_0^y
    \frac{\dot\sigma_\theta^v(r)}{\sigma_\theta(r)^2}\dd r,
  \label{eq:parametric-F-derivative}\\
  \dot G_\theta^v(z)
  &=-\sigma_\theta(G_\theta(z))
    \dot F_\theta^v(G_\theta(z)).
  \label{eq:parametric-G-derivative}
\end{align}
For \(g=G_\theta(z)\) and \(\dot g=\dot G_\theta^v(z)\), the derivatives of
the transformed coefficients are
\begin{align}
  \dot{\widehat b}_\theta^v(z)
  ={}&
  \frac{\dot b_\theta^v(g)+b_\theta'(g)\dot g}{\sigma_\theta(g)}
  -\frac{b_\theta(g)}{\sigma_\theta(g)^2}
  \bigl(\dot\sigma_\theta^v(g)+\sigma_\theta'(g)\dot g\bigr),
  \label{eq:parametric-bhat-derivative}\\
  \dot{\widehat c}_\theta^v(z)
  ={}&-\frac12
  \bigl(\dot\sigma_\theta^{v\prime}(g)+\sigma_\theta''(g)\dot g\bigr).
  \label{eq:parametric-chat-derivative}
\end{align}

\begin{lemma}[Pathwise parameter variation]
\label{lem:parametric-variation}
Under \Cref{ass:coefficients,ass:parametric}, fix \(x\in\Omega\),
\(q\in\cA_\Lambda\), \(\theta\in\operatorname{int}\Theta_1\),
\(s\in[0,T]\), and \(y\in\R\).
Set \(Z=Z^{\theta,s,y}(x,q)\).  For every \(v\in\R^d\), the derivative
\(W^v=D_\theta Z^{\theta,s,y}(x,q)[v]\) exists in \(C([s,T])\) and is the
unique solution of
\begin{align}
  W_t^v
  ={}&\dot F_\theta^v(y)
  +\int_s^t
    \left(
      \widehat b_\theta'(Z_r)W_r^v
      +\dot{\widehat b}_\theta^v(Z_r)
    \right)\dd r
  \notag\\
  &+\int_s^t
    \left(
      \widehat c_\theta'(Z_r)W_r^v
      +\dot{\widehat c}_\theta^v(Z_r)
    \right)\dd q_r.
  \label{eq:parametric-Lamperti-variation}
\end{align}
If
\[
  \mathcal J^Z_{s,t}
  =\exp\left\{
    \int_s^t\widehat b_\theta'(Z_r)\dd r
    +\int_s^t\widehat c_\theta'(Z_r)\dd q_r
  \right\},
\]
then
\begin{align}
  W_t^v
  =\mathcal J^Z_{s,t}\bigg[
    \dot F_\theta^v(y)
    &+\int_s^t(\mathcal J^Z_{s,r})^{-1}
      \dot{\widehat b}_\theta^v(Z_r)\dd r
    \notag\\
    &+\int_s^t(\mathcal J^Z_{s,r})^{-1}
      \dot{\widehat c}_\theta^v(Z_r)\dd q_r
  \bigg].
  \label{eq:parametric-variation-of-constants}
\end{align}
Writing \(Y_t=\Phi_{s,t}^\theta(x,q;y)\), the derivative of the original
flow is
\begin{equation}\label{eq:parametric-flow-derivative}
  D_\theta\Phi_{s,t}^\theta(x,q;y)[v]
  =\dot G_\theta^v(Z_t)+G_\theta'(Z_t)W_t^v
  =\sigma_\theta(Y_t)\bigl(W_t^v-\dot F_\theta^v(Y_t)\bigr).
\end{equation}
The derivative is jointly continuous in all its arguments.  For every compact
\(\Theta_0\Subset\operatorname{int}\Theta_1\) and \(m<\infty\),
\begin{equation}\label{eq:parametric-pathwise-bound}
  \sup_{\substack{q\in\cA_\Lambda,\ \theta\in\Theta_0,\
  0\le s\le t\le T\\
  |y|\le m}}
  \abs{D_\theta\Phi_{s,t}^\theta(x,q;y)[v]}
  \le C_{m,\Theta_0}(1+\norm{x}_\infty)\abs{v},
  \qquad x\in\Omega.
\end{equation}
\end{lemma}

\begin{proof}
The Banach-valued \(C^1\) assumption, the identity
\(F_\theta(G_\theta(z))=z\), and the chain rule give
\eqref{eq:parametric-F-derivative}--\eqref{eq:parametric-chat-derivative}.
Subtracting the finite-measure equations at \(\theta+\eps v\) and \(\theta\),
dividing by \(\eps\), and applying the mean-value formula gives
\eqref{eq:parametric-Lamperti-variation} by Gronwall's inequality for
\(\dd r+\dd q_r\).  Variation of constants gives
\eqref{eq:parametric-variation-of-constants}, and the chain rule gives
\eqref{eq:parametric-flow-derivative}.

The transformed solution obeys
\[
  \sup_{s\le r\le T}\abs{Z_r}
  \le C(1+\abs{y}+\norm{x}_\infty),
\]
uniformly in \(q\in\cA_\Lambda\).  Equations
\eqref{eq:parametric-bhat-derivative}--\eqref{eq:parametric-chat-derivative}
have at most linear growth in \(Z\), while the undotted spatial derivatives
are bounded.  Formula \eqref{eq:parametric-variation-of-constants} proves
\eqref{eq:parametric-pathwise-bound}.  Joint continuity follows from the
same Gronwall and Helly--Bray argument as in
\Cref{lem:deterministic-flow}, applied simultaneously to the equation for
\(Z\) and its linear variation.
\end{proof}

\begin{theorem}[Common coefficient-parameter tangent]
\label{thm:parametric-tangent}
Under \Cref{ass:coefficients,ass:parametric}, define
\[
  D_\theta S_{s,t}^\theta(x;y)[v]
  =D_\theta\Phi_{s,t}^\theta(x,Q(x);y)[v].
\]
Then this map is Borel, causal, linear in \(v\), and continuous in
\((\theta,s,t,y,v)\) for every fixed \(x\).  It satisfies the affine tangent
cocycle
\begin{align}
  D_\theta S_{s,t}^\theta(x;y)[v]
  ={}&D_\theta S_{u,t}^\theta
  \bigl(x;S_{s,u}^\theta(x;y)\bigr)[v]
  \notag\\
  &+\left.
    \partial_\xi S_{u,t}^\theta(x;\xi)
    \right|_{\xi=S_{s,u}^\theta(x;y)}
    D_\theta S_{s,u}^\theta(x;y)[v],
  \qquad s\le u\le t.
  \label{eq:parametric-affine-cocycle}
\end{align}

For the common raw flow, define
\begin{equation}\label{eq:parametric-flow-field}
  \mathscr F_m^{\mathrm{par}}(x)
  =
  \left(
    S_{s,t}^\theta(x;y),
    \partial_yS_{s,t}^\theta(x;y),
    \nabla_\theta S_{s,t}^\theta(x;y)
  \right)_{\theta\in\Theta,\,s\le t,\,|y|\le m}.
\end{equation}
Then \(\mathscr F_m^{\mathrm{par}}\) is Borel with values in
\begin{equation}\label{eq:parametric-field-space}
  E_m^{\mathrm{par}}
  =C\bigl(
    \Theta\times\{(s,t):s\le t\}\times[-m,m];\R^{2+d}
  \bigr)
\end{equation}
and is continuous on every compact capacity core \(K_R\) from
\Cref{thm:causal-bracket}.

Finally, under each \(P\in\mathfrak M_\Lambda\), the process
\[
  \Xi_t^{\theta,v}
  :=D_\theta S_{s,t}^\theta(X;y)[v]
\]
is the unique adapted solution of the stochastic variational equation
\begin{align}
  \Xi_t^{\theta,v}
  ={}&\int_s^t
    \left(
      b_\theta'(Y_r)\Xi_r^{\theta,v}
      +\dot b_\theta^v(Y_r)
    \right)\dd r
  \notag\\
  &+\int_s^t
    \left(
      \sigma_\theta'(Y_r)\Xi_r^{\theta,v}
      +\dot\sigma_\theta^v(Y_r)
    \right)\dd X_r,
  \qquad
  Y_r=S_{s,r}^\theta(X;y).
  \label{eq:modelwise-parametric-variation}
\end{align}
It is adapted and has finite moments of every order.  More precisely, for
every \(p\ge2\) and \(m<\infty\),
\begin{equation}
  \sup_{P\in\mathfrak M_\Lambda}
  E^P\left[
    \sup_{\substack{\theta\in\Theta,\ 0\le s\le t\le T\\ |y|\le m}}
    \norm{\nabla_\theta S_{s,t}^\theta(X;y)}^p
  \right]<\infty,
  \label{eq:uniform-parametric-moments}
\end{equation}
\end{theorem}

\begin{proof}
Composition with the Borel causal map \(x\mapsto(x,Q(x))\) and
\Cref{lem:parametric-variation} give the first assertions.  Differentiating
the flow cocycle in \(\theta\), including the parameter-dependent
intermediate state, gives \eqref{eq:parametric-affine-cocycle}.

For a fixed core, the set
\(\{(x,Q(x)):x\in K_R\}\) is compact.  Joint continuity in
\Cref{lem:parametric-variation} proves continuity of
\(\mathscr F_m^{\mathrm{par}}|_{K_R}\), and joint Borel measurability plus
continuity in the field variables gives Borel measurability into the
separable space \(E_m^{\mathrm{par}}\).

Fix \(P\in\mathfrak M_\Lambda\).  For each \(t\), causality gives the
stopped-path identity
\[
  D_\theta S_{s,t}^\theta(x;y)[v]
  =D_\theta S_{s,t}^\theta(x_{\cdot\wedge t};y)[v].
\]
Together with Borel measurability this makes the random variable on the
left measurable with respect to the raw canonical \(\sigma\)-field at time
\(t\).  Its continuous paths are therefore progressively measurable under
the completed, right-continuous filtration of every model.

Fix \((\theta,s,y,v)\) and let \(\eps\to0\) along values for which
\(\theta+\lambda\eps v\in\Theta_1\) for every \(0\le\lambda\le1\).  Under
\(P\), the common flows at
\(\theta+\eps v\) and \(\theta\) are their classical strong solutions.
Write these solutions as \(Y^\eps\) and \(Y\), and put
\(\Delta^\eps=(Y^\eps-Y)/\eps\).  The mean-value formula gives
\[
  \Delta_t^\eps
  =\int_s^t(A_r^\eps\Delta_r^\eps+\beta_r^\eps)\dd r
   +\int_s^t(C_r^\eps\Delta_r^\eps+\gamma_r^\eps)\dd X_r,
\]
where
\begin{align*}
A_r^\eps
&=\int_0^1 b_{\theta+\eps v}'
  \bigl(Y_r+\lambda(Y_r^\eps-Y_r)\bigr)\dd\lambda,
&
\beta_r^\eps
&=\frac{b_{\theta+\eps v}(Y_r)-b_\theta(Y_r)}{\eps},
\\
C_r^\eps
&=\int_0^1 \sigma_{\theta+\eps v}'
  \bigl(Y_r+\lambda(Y_r^\eps-Y_r)\bigr)\dd\lambda,
&
\gamma_r^\eps
&=\frac{\sigma_{\theta+\eps v}(Y_r)-\sigma_\theta(Y_r)}{\eps}.
\end{align*}
The uniform coefficient bounds and BDG--Gronwall imply, for \(p\ge2\),
\[
  \sup_{0<|\eps|\le\eps_0}
  E^P\!\left[\sup_{s\le r\le T}|\Delta_r^\eps|^p\right]<\infty,
  \qquad
  E^P\!\left[\sup_{s\le r\le T}|Y_r^\eps-Y_r|^p\right]
  =O(|\eps|^p).
\]
The Banach-valued \(C^1\) assumption and the bounded second spatial
derivatives yield convergence of
\[
  (A^\eps,\beta^\eps,C^\eps,\gamma^\eps)
  \quad\hbox{to}\quad
  (b_\theta'(Y),\dot b_\theta^v(Y),
   \sigma_\theta'(Y),\dot\sigma_\theta^v(Y))
\]
in the norms entering the drift and martingale integrals.  A second
BDG estimate, using
\(\dd\langle X\rangle_r^P\le\Lambda\dd r\), gives for the solution
\(\Xi\) of \eqref{eq:modelwise-parametric-variation}
\[
 E^P\!\left[\sup_{s\le u\le t}|\Delta_u^\eps-\Xi_u|^p\right]
 \le C_p\int_s^t
 E^P\!\left[\sup_{s\le u\le r}|\Delta_u^\eps-\Xi_u|^p\right]\dd r
 +\rho_p(\eps),
 \qquad \rho_p(\eps)\longrightarrow0.
\]
Gronwall's inequality shows that \(\Delta^\eps\) converges in
\(L^p(P;C([s,T]))\) to the unique solution of
\eqref{eq:modelwise-parametric-variation}.  The same quotients converge
pathwise by \Cref{lem:parametric-variation}; uniqueness of limits in
probability identifies the two limits.  Finally,
\eqref{eq:parametric-pathwise-bound} with a compact neighborhood of
\(\Theta\), followed by the uniform moment bound for \(\norm X_\infty\),
proves \eqref{eq:uniform-parametric-moments}.  Thus the supremum over the
field parameters in \eqref{eq:uniform-parametric-moments} is controlled
before expectation is taken.
\end{proof}

In this parameter calculus, \(\mathfrak M_\Lambda\) and the bracket map
\(Q\) are fixed; the parameter enters through \((b_\theta,\sigma_\theta)\).
Consequently, the tangent equation has no bracket-variation term.

For \(\eta=(v,w)\in\R^d\times\R\), whenever
\(\theta+[0,1]v\subset\operatorname{int}\Theta_1\), define the joint tangent
and remainder by
\begin{align}
  \mathcal T_{s,t}^{\theta,y}(x;\eta)
  &:=D_\theta S_{s,t}^\theta(x;y)[v]
    +w\,\partial_yS_{s,t}^\theta(x;y),
  \label{eq:joint-flow-tangent}\\
  \mathcal R_{s,t}^{\theta,y}(x;\eta)
  &:=S_{s,t}^{\theta+v}(x;y+w)-S_{s,t}^\theta(x;y)
    -\mathcal T_{s,t}^{\theta,y}(x;\eta).
  \label{eq:joint-flow-remainder}
\end{align}

\begin{proposition}[Uniform joint first-order expansion]
\label{prop:uniform-joint-expansion}
Under \Cref{ass:parametric}, let
\(\Theta_0\Subset\operatorname{int}\Theta_1\) be compact and choose
\(\rho>0\) such that
\[
  (\Theta_0)_\rho
  :=\{\vartheta:\operatorname{dist}(\vartheta,\Theta_0)\le\rho\}
  \subset\operatorname{int}\Theta_1.
\]
For every \(R,m<\infty\),
\begin{equation}\label{eq:core-uniform-Frechet}
  \lim_{\delta\downarrow0}
  \sup_{\substack{x\in K_R,\ \theta\in\Theta_0,\ |y|\le m,\ s\le t\\
  0<\abs{\eta}\le\min\{\delta,\rho,1\}}}
  \frac{\abs{\mathcal R_{s,t}^{\theta,y}(x;\eta)}}{\abs{\eta}}=0.
\end{equation}
For every \(1\le p<\infty\), the corresponding expansion is uniform over
the nondominated model class:
\begin{equation}\label{eq:robust-uniform-Frechet}
  \lim_{\delta\downarrow0}
  \sup_{\substack{P\in\mathfrak M_\Lambda,\ \theta\in\Theta_0,\ |y|\le m,\
  0\le s\le T\\0<\abs{\eta}\le\min\{\delta,\rho,1\}}}
  \frac{1}{\abs{\eta}}
  \left(
    E^P\left[
      \sup_{s\le t\le T}
      \abs{\mathcal R_{s,t}^{\theta,y}(X;\eta)}^p
    \right]
  \right)^{1/p}=0.
\end{equation}
\end{proposition}

\begin{proof}
Let
\[
  \Delta_T=\{(s,t):0\le s\le t\le T\},
  \qquad \mathbb B=C(\Delta_T),
\]
and define the Banach-valued map
\[
  F(x,\theta,y)=\bigl(S_{s,t}^\theta(x;y)\bigr)_{(s,t)\in\Delta_T}.
\]
For fixed \(x\), joint continuity of the pointwise derivative in
\((\theta,y,s,t)\), compactness of \(\Delta_T\), and the mean-value formula
give
\[
  \frac{\norm{F(x,z+h)-F(x,z)-D_zF(x,z)h}_{\mathbb B}}{|h|}
  \le
  \sup_{\lambda\in[0,1]}
  \norm{D_zF(x,z+\lambda h)-D_zF(x,z)}
  \longrightarrow0.
\]
Thus \(z\mapsto F(x,z)\) is \(C^1\) as a \(\mathbb B\)-valued map, and
its derivative is precisely the field
\(\mathcal T^{\theta,y}(x;\cdot)\).  Set
\[
  C_m=\Theta_0\times[-m,m],
  \qquad \rho_*=\min\{\rho,1\}.
\]
Then
\[
  (C_m)_{\rho_*}
  \subset(\Theta_0)_\rho\times[-m-1,m+1]
  \subset\operatorname{int}\Theta_1\times\R.
\]
By
\Cref{thm:common-flow,thm:parametric-tangent}, the first jet of \(F\) is
jointly continuous on every
\[
  K_R\times(\Theta_0)_\rho\times[-m-1,m+1].
\]
Moreover, \eqref{eq:uniform-flow-Jacobian} and
\eqref{eq:parametric-pathwise-bound} give the first-jet envelope
\begin{equation}\label{eq:joint-first-jet-envelope}
  \sup_{\substack{\vartheta\in(\Theta_0)_\rho,\ |z|\le m+1}}
  \norm{D_{(\vartheta,z)}F(x,\vartheta,z)}
  \le C_{m,\Theta_0}(1+\norm{x}_\infty).
\end{equation}
The bracket bound and the Burkholder--Davis--Gundy inequality make the
right-hand side uniformly integrable to every finite power over
\(P\in\mathfrak M_\Lambda\).  The core modulus in
\Cref{thm:core-C1-transfer} gives
\eqref{eq:core-uniform-Frechet}; applying
\eqref{eq:C1-transfer-bound} with an arbitrary \(q>p\) gives
\eqref{eq:robust-uniform-Frechet}.
\end{proof}

\section{Law-level stability and robust sensitivity}
\label{sec:probability-output}

The common compact cores now turn the deterministic flow structure into
statements about the entire nondominated class.

\begin{theorem}[Wasserstein stability of the flow fields]
\label{thm:flow-law-transfer}
Fix \(m<\infty\).  For every \(1\le p<\infty\), if
\(P_n,P\in\mathfrak M_\Lambda\) and \(P_n\Rightarrow P\), then
\begin{equation}\label{eq:flow-field-Wp}
  W_p\bigl(P_n\circ\mathscr F_m^{-1},
            P\circ\mathscr F_m^{-1}\bigr)\longrightarrow0.
\end{equation}
Consequently,
\[
  \{P\circ\mathscr F_m^{-1}:P\in\mathfrak M_\Lambda\}
\]
is compact in \(\mathfrak P_p(E_m)\) endowed with \(W_p\).

Under \Cref{ass:parametric}, the same assertions hold with
\(\mathscr F_m^{\mathrm{par}}\) and \(E_m^{\mathrm{par}}\) in place of
\(\mathscr F_m\) and \(E_m\).  For either field \(\mathscr G_m\) and
every \(1\le r<\infty\),
\begin{equation}\label{eq:parametric-field-moments}
  \sup_{P\in\mathfrak M_\Lambda}
  E^P\bigl[\norm{\mathscr G_m}^{r}\bigr]<\infty.
\end{equation}
Hence, if \(\Psi\) is continuous on the corresponding field space and
\[
  |\Psi(e)|\le C(1+\norm e^r)
\]
for some \(r\ge0\), then there exists \(P^*\in\mathfrak M_\Lambda\)
such that
\begin{equation}\label{eq:field-worst-case}
  \sup_{P\in\mathfrak M_\Lambda}E^P[\Psi(\mathscr G_m)]
  =E^{P^*}[\Psi(\mathscr G_m)].
\end{equation}
More generally, all parametric assertions remain valid when the field in
\eqref{eq:parametric-flow-field} is indexed by any compact set
\(\Theta_0\Subset\operatorname{int}\Theta_1\) in place of \(\Theta\).
\end{theorem}

\begin{proof}
Core continuity follows from
\Cref{thm:common-flow,thm:parametric-tangent}.  The deterministic flow
estimates, \eqref{eq:uniform-flow-Jacobian}, and
\eqref{eq:parametric-pathwise-bound} imply, for the relevant field,
\[
  \norm{\mathscr G_m}\le C_m(1+\norm X_\infty).
\]
The Burkholder--Davis--Gundy inequality and
\(\langle X\rangle_T\le\Lambda T\) prove
\eqref{eq:parametric-field-moments} for every finite \(r\).

By \Cref{thm:core-transfer}(ii), weak convergence of the input laws gives
weak convergence of the field laws.  Uniform integrability of their
\(p\)-th moments, obtained by choosing any exponent larger than \(p\),
upgrades this to \(W_p\)-convergence.  Weak compactness of
\(\mathfrak M_\Lambda\) from \Cref{lem:model-compactness} then makes each
image family \(W_p\)-compact.  Finally, if \(\Psi\) has polynomial growth,
choose a moment order larger than its growth order.  Truncation and uniform
integrability show that
\(P\mapsto E^P[\Psi(\mathscr G_m)]\) is continuous on
\(\mathfrak M_\Lambda\), and compactness gives \eqref{eq:field-worst-case}.
\end{proof}

\subsection{Robust sensitivity}

The next proposition is the compact-index Berge--Danskin theorem used
below; see \cite{BonnansShapiro00}.

\begin{proposition}[Compact Berge--Danskin principle]
\label{thm:compact-robust-envelope}
Let \(\mathcal M\) be a compact metric space, let \(U\subset\R^d\) be open,
and let \(g:\mathcal M\times U\to\R\) be continuous.  Suppose that
\(g(m,\cdot)\) is \(C^1\) and \((m,z)\mapsto D_zg(m,z)\) is continuous.
Set
\[
  V(z)=\max_{m\in\mathcal M}g(m,z),
  \qquad \mathcal A(z)=\Argmax_{m\in\mathcal M}g(m,z).
\]
Then \(V\) is continuous, \(\mathcal A\) has nonempty compact values,
closed graph, and is upper hemicontinuous.  Moreover,
\begin{equation}\label{eq:abstract-Danskin}
  V'_{H,+}(z;h)
  :=\lim_{\substack{t_n\downarrow0,\ h_n\to h\\z+t_nh_n\in U}}
  \frac{V(z+t_nh_n)-V(z)}{t_n}
  =\max_{m\in\mathcal A(z)}D_zg(m,z)h.
\end{equation}
If all active gradients equal \(\ell\in(\R^d)^*\), then \(V\) is
Fr\'echet differentiable at \(z\) with \(D V(z)=\ell\).  If
\(\mathcal A(z)=\{m^*\}\), then \(z_n\to z\) and
\(m_n\in\mathcal A(z_n)\) imply \(m_n\to m^*\).
\end{proposition}

\begin{proof}
On compact subsets of \(U\), joint continuity of \(D_zg\) gives a
first-order remainder uniform in \(m\).  Compactness gives continuity of
\(V\) and the stated properties of \(\mathcal A\).  Testing with a fixed
active \(m\) yields the lower bound in \eqref{eq:abstract-Danskin}.  For
the upper bound, choose \(m_n\in\mathcal A(z+t_nh_n)\) and pass to an
active limit point.  The same argument shows that active gradients near
\(z\) converge uniformly to \(\ell\) when they agree at \(z\), giving the
Fr\'echet expansion; uniqueness identifies every cluster point of
\((m_n)\).
\end{proof}

The extension in \Cref{ass:parametric} defines the flow and its payoff on
\(\operatorname{int}\Theta_1\).  Fix \((s,t)\) and \(f\in C^1(\R)\) such
that, for some \(k\ge0\),
\begin{equation}\label{eq:payoff-polynomial-growth}
  |f(r)|+|f'(r)|\le C(1+|r|^k),
  \qquad r\in\R.
\end{equation}
Write
\(z=(\theta,y)\in\operatorname{int}\Theta_1\times\R\), and set
\[
  g(P,z)=E^P[f(S_{s,t}^\theta(X;y))],
  \qquad
  V(z)=\max_{P\in\mathfrak M_\Lambda}g(P,z),
  \qquad
  \cP^*(z)=\Argmax_{P\in\mathfrak M_\Lambda}g(P,z).
\]

\begin{theorem}[Joint parametric robust sensitivity]
\label{thm:parametric-robust-sensitivity}
Assume \Cref{ass:parametric}.  The value function \(V\) is continuous.  The
correspondence \(z\mapsto\cP^*(z)\) has nonempty weakly compact values,
closed graph, and is upper hemicontinuous.  The map \(g\) is jointly
continuous in \((P,z)\), continuously differentiable in \(z\), and
\begin{align}
  D_zg(P,z)[(v,w)]
  =E^P\Bigl[
    f'(S_{s,t}^\theta(X;y))
    \bigl(
      D_\theta S_{s,t}^\theta(X;y)[v]
      +w\,\partial_yS_{s,t}^\theta(X;y)
    \bigr)
  \Bigr].
  \label{eq:joint-modelwise-sensitivity}
\end{align}
For every compact
\(C\subset\operatorname{int}\Theta_1\times\R\), choose \(\rho>0\) such that
the closed \(\rho\)-neighborhood \(C_\rho\) is contained in
\(\operatorname{int}\Theta_1\times\R\).  The first-order expansion is
uniform:
\begin{equation}\label{eq:uniform-payoff-Frechet}
  \lim_{\delta\downarrow0}
  \sup_{\substack{P\in\mathfrak M_\Lambda,\ z\in C,\ z+\eta\in C_\rho\\
  0<\abs{\eta}\le\delta}}
  \frac{\abs{g(P,z+\eta)-g(P,z)-D_zg(P,z)[\eta]}}{\abs{\eta}}=0.
\end{equation}

For every \(\eta=(v,w)\in\R^{d+1}\), \(V\) has the one-sided Hadamard
directional derivative
\begin{align}
  V'_{H,+}(z;\eta)
  &:={}
  \lim_{\substack{\eps\downarrow0,\ \eta'\to\eta\\
  z+\eps\eta'\in\operatorname{int}\Theta_1\times\R}}
  \frac{V(z+\eps\eta')-V(z)}{\eps}
  \notag\\
  &=\max_{P\in\cP^*(z)}
  E^P\Bigl[
    f'(S_{s,t}^\theta(X;y))
    \bigl(
      D_\theta S_{s,t}^\theta(X;y)[v]
      +w\,\partial_yS_{s,t}^\theta(X;y)
    \bigr)
  \Bigr].
  \label{eq:joint-Danskin}
\end{align}
If all \(P\in\cP^*(z)\) give the same vector
\begin{equation}\label{eq:common-active-gradient}
  E^P\left[
    f'(S_{s,t}^\theta(X;y))
    \bigl(
      \nabla_\theta S_{s,t}^\theta(X;y),
      \partial_yS_{s,t}^\theta(X;y)
    \bigr)
  \right]\in\R^{d+1},
\end{equation}
then \(V\) is Fr\'echet differentiable at \(z\), with derivative equal to
this common vector.  This condition holds in particular when the maximizing
law is unique.  If \(z_n\to z\), \(P_n\in\cP^*(z_n)\), and
\(\cP^*(z)=\{P^*\}\), then \(P_n\Rightarrow P^*\).
\end{theorem}

\begin{proof}
Fix a compact \(C\subset\operatorname{int}\Theta_1\times\R\).  Choose
\(\Theta_0\Subset\operatorname{int}\Theta_1\) and \(m<\infty\) so that
\(\Theta_0\times[-m,m]\) contains a neighborhood of \(C\).  The local
parametric assertion in \Cref{thm:flow-law-transfer} then applies to the
full first-jet field on this enlarged index set.

For \(z=(\theta,y)\) and \(\eta=(v,w)\), abbreviate
\[
  Y=S_{s,t}^\theta(X;y),\qquad
  \Delta_\eta=S_{s,t}^{\theta+v}(X;y+w)-Y,\qquad
  A\eta=\mathcal T_{s,t}^{\theta,y}(X;\eta),\qquad
  R_\eta=\Delta_\eta-A\eta.
\]
Taylor's formula gives
\begin{align*}
&\frac{|f(Y+\Delta_\eta)-f(Y)-f'(Y)A\eta|}{|\eta|}\\
&\quad\le
 |f'(Y)|\frac{|R_\eta|}{|\eta|}
 +\frac{|\Delta_\eta|}{|\eta|}
  \int_0^1|f'(Y+r\Delta_\eta)-f'(Y)|\,\dd r.
\end{align*}
The first term tends to zero in expectation, uniformly in
\(P\in\mathfrak M_\Lambda\) and \(z\in C\), by H\"older's inequality,
\eqref{eq:robust-uniform-Frechet}, and the moments of every order.  For the
second term, \(|\Delta_\eta|/|\eta|\) is bounded in every finite
\(L^p(P)\), uniformly in \(P\) and \(z\in C\).  On bounded state sets use
uniform continuity of \(f'\); on their complements use
\[
  |f'(a)-f'(b)|\le C(2+|a|^k+|b|^k)
\]
and the higher moment bounds.  This proves
\eqref{eq:joint-modelwise-sensitivity} and
\eqref{eq:uniform-payoff-Frechet}.

Choose a Wasserstein order strictly larger than \(k+1\).  The integrands in
\(g\) and \(D_zg\) are continuous functions of the local first-jet field
with growth of order at most \(k+1\).  The local form of
\Cref{thm:flow-law-transfer} and uniform integrability therefore show that
both \((P,z)\mapsto g(P,z)\) and
\((P,z)\mapsto D_zg(P,z)\) are jointly continuous.

The space \(\mathfrak M_\Lambda\) is compact by
\Cref{lem:model-compactness}.  Apply
\Cref{thm:compact-robust-envelope} with
\(\mathcal M=\mathfrak M_\Lambda\).  Its four conclusions give,
respectively, continuity and stability of the active laws,
\eqref{eq:joint-Danskin}, the common-gradient Fr\'echet criterion, and the
convergence of maximizers under uniqueness.
Since \(\Omega\) is Polish, the weak topology on
\(\mathfrak P(\Omega)\) is metrizable, so the last sequential statement is
equivalent to weak convergence of the whole sequence.
\end{proof}

\begin{example}[Active-model switching]\label{ex:active-model-switch}
The preceding proof applies to any weakly compact subfamily of
\(\mathfrak M_\Lambda\).  Assume \(T>0\) and take \(d=1\),
\[
  b_\theta(r)=\theta,\qquad \sigma_\theta(r)=1,
\]
with \(\theta\) in a compact interval contained in an open parameter
interval, and set \(s=0\), \(t=T\).  Then
\[
  S_{0,T}^\theta(X;y)=m+X_T,
  \qquad m=y+\theta T.
\]
For \(0\le a\le\Lambda\), let
\(P_a=\operatorname{Law}(\sqrt a\,W)\) on \(\Omega\), where \(W\) is a
standard Brownian motion, and consider
\[
  \mathfrak M_\Lambda^{\mathrm{cv}}
  =\{P_a:0\le a\le\Lambda\}.
\]
The map \(a\mapsto P_a\) is weakly continuous, so this is a weakly compact
subfamily of \(\mathfrak M_\Lambda\).  Distinct members are mutually
singular because \(Q_t(X)=at\) for all \(t\in[0,T]\),
\(P_a\)-almost surely.  Choose
\[
  f(r)=r^4-cr^2,\qquad c>3\Lambda T,
\]
and write \(\tau=\Lambda T\).  Since
\(X_T\sim N(0,aT)\) under \(P_a\),
\[
  g(P_a,\theta,y)
  =m^4-cm^2+(6m^2-c)aT+3(aT)^2.
\]
The right-hand side is strictly convex in \(aT\), so its maximum over
\(0\le a\le\Lambda\) is attained at \(a=0\) or \(a=\Lambda\).  Moreover,
\[
  g(P_\Lambda,\theta,y)-g(P_0,\theta,y)
  =\tau(6m^2+3\tau-c).
\]
Thus, for
\[
  \rho=\sqrt{\frac{c-3\tau}{6}},
\]
the active set over \(\mathfrak M_\Lambda^{\mathrm{cv}}\) is
\[
  \cP_{\mathrm{cv}}^*(\theta,y)
  =
  \begin{cases}
    \{P_0\},& |m|<\rho,\\
    \{P_0,P_\Lambda\},& |m|=\rho,\\
    \{P_\Lambda\},& |m|>\rho.
  \end{cases}
\]
At either switching surface \(m=\pm\rho\), the active gradients differ by
\[
  \nabla_{(\theta,y)}g(P_\Lambda,\theta,y)
  -\nabla_{(\theta,y)}g(P_0,\theta,y)
  =12\tau m\,(T,1)\ne0.
\]
The robust value \(V_{\mathrm{cv}}=\max_{0\le a\le\Lambda}g(P_a,\cdot)\)
is therefore not Fr\'echet differentiable on the switching surfaces.  If
\(d_\eta=Tv+w\), its one-sided Hadamard derivative there is
\[
  V'_{\mathrm{cv},H,+}((\theta,y);(v,w))
  =
  \max\left\{
    (4m^3-2cm)d_\eta,\,
    \bigl(4m^3+2(6\tau-c)m\bigr)d_\eta
  \right\},
\]
in agreement with \eqref{eq:joint-Danskin}.
\end{example}

\section{Finite-variation tangents and Gaussian identification}
\label{sec:tangent}

Finite-variation invariance of the bracket yields a pathwise tangent in the
driving signal.

\subsection{Finite-variation response}

Let \(C_0^{\mathrm{bv}}([0,T])\) be the continuous finite-variation paths
starting at zero.

\begin{proposition}[Finite-variation invariance]\label{prop:bv-invariance}
For every \(x\in\Omega\) and every
\(h\in C_0^{\mathrm{bv}}([0,T])\),
\begin{equation}\label{eq:Q-bv-invariance}
  Q(x+h)=Q(x).
\end{equation}
\end{proposition}

\begin{proof}
Let \(\omega_x(\delta)\) and \(\omega_h(\delta)\) denote the uniform
moduli of continuity.  Expanding the completed sums gives, uniformly in
\(t\),
\begin{align*}
  \abs{V_t^n(x+h)-V_t^n(x)}
  &\le
  2\max_k\abs{x_{t_{k+1}^n}-x_{t_k^n}}
    \sum_k\abs{h_{t_{k+1}^n}-h_{t_k^n}}\\
  &\quad+
  \max_k\abs{h_{t_{k+1}^n}-h_{t_k^n}}
    \sum_k\abs{h_{t_{k+1}^n}-h_{t_k^n}}\\
  &\le
  \bigl(2\omega_x(T2^{-n})+\omega_h(T2^{-n})\bigr)\Var(h).
\end{align*}
The right-hand side tends to zero.  Since
\(\mathcal R_\Lambda\) is nonexpansive,
\[
  \norm{q^n(x+h)-q^n(x)}_\infty\longrightarrow0.
\]
The two pointwise limsup functions in \eqref{eq:total-bracket} are therefore
equal.
\end{proof}

Thus finite-variation directions change the first-order signal but do not
change its second-order coordinate.

Fix \((x,\theta,s,y)\), write
\[
  q=Q(x),
  \qquad
  Z_t=Z_t^{\theta,s,y}(x,q),
  \qquad
  Y_t=S_{s,t}^\theta(x;y).
\]

\begin{theorem}[Common finite-variation directional derivative]
\label{thm:bv-derivative}
For every \(h\in C_0^{\mathrm{bv}}([0,T])\), the map
\[
  \eps\longmapsto S_{s,t}^\theta(x+\eps h;y)
\]
is differentiable at \(\eps=0\), locally uniformly in
\((\theta,s,t,y)\).  In Lamperti coordinates, \(U=D_hZ\) is the unique
solution of
\begin{equation}\label{eq:Lamperti-directional-equation}
  U_t
  =h_t-h_s
  +\int_s^t\widehat b_\theta'(Z_r)U_r\dd r
  +\int_s^t\widehat c_\theta'(Z_r)U_r\dd Q_r(x).
\end{equation}
Moreover,
\begin{equation}\label{eq:S-directional-from-U}
  D_hS_{s,t}^\theta(x;y)
  =G_\theta'(Z_t)U_t.
\end{equation}

Define the pathwise response kernel
\begin{equation}\label{eq:response-kernel}
  \mathcal K_{s,t}^\theta(r;x,y)
  =
  \sigma_\theta(Y_r)
  \left.
    \partial_\xi S_{r,t}^\theta(x;\xi)
  \right|_{\xi=Y_r},
  \qquad s\le r\le t.
\end{equation}
Then
\begin{equation}\label{eq:response-representation}
  D_hS_{s,t}^\theta(x;y)
  =
  \int_{(s,t]}
  \mathcal K_{s,t}^\theta(r;x,y)\dd h_r,
\end{equation}
and
\begin{equation}\label{eq:response-kernel-bound}
  \sup_{\substack{x\in\Omega,\ \theta\in\Theta,\ y\in\R\\
                   0\le s\le r\le t\le T}}
  \abs{\mathcal K_{s,t}^\theta(r;x,y)}
  \le\overline\sigma C_*.
\end{equation}
Consequently,
\[
  \abs{D_hS_{s,t}^\theta(x;y)}
  \le\overline\sigma C_*\Var(h;[s,t]).
\]
\end{theorem}

\begin{proof}
By \Cref{prop:bv-invariance},
\(Q(x+\eps h)=Q(x)\) for every \(\eps\).  Write \(Z^\eps\) for the
Lamperti solution driven by \(x+\eps h\) and set
\(U^\eps=(Z^\eps-Z)/\eps\).  The mean-value formula gives
\begin{align*}
  U_t^\eps
  ={}&h_t-h_s
  +\int_s^t B_r^\eps U_r^\eps\dd r
  +\int_s^t C_r^\eps U_r^\eps\dd Q_r(x),
\end{align*}
where
\[
  B_r^\eps
  =\int_0^1\widehat b_\theta'
    \bigl(Z_r+\lambda(Z_r^\eps-Z_r)\bigr)\dd\lambda,
  \qquad
  C_r^\eps
  =\int_0^1\widehat c_\theta'
    \bigl(Z_r+\lambda(Z_r^\eps-Z_r)\bigr)\dd\lambda.
\]
The common derivative bounds and Gronwall's inequality for
\(\dd r+\dd Q_r\) give a uniform bound on \(U^\eps\) and first imply
\(\sup_t\abs{Z_t^\eps-Z_t}=O(\abs\eps)\).  Since
\(\widehat b_\theta''\) and \(\widehat c_\theta''\) are uniformly bounded,
\[
  \sup_{\theta,\,0\le s\le t\le T,\,|y|\le m}
  \abs{U_t^\eps-U_t}
  \longrightarrow0
\]
for every \(m<\infty\), after a second Gronwall estimate.  This proves local
uniform differentiability.  Applying the same argument to \(G_\theta\),
whose second derivative is uniformly bounded, gives
\eqref{eq:S-directional-from-U}.

Variation of constants in the scalar linear equation yields
\[
  U_t
  =\int_{(s,t]}
  \exp\left\{
    \int_r^t\widehat b_\theta'(Z_u)\dd u
    +\int_r^t\widehat c_\theta'(Z_u)\dd Q_u(x)
  \right\}\dd h_r.
\]
Multiplication by
\(G_\theta'(Z_t)=\sigma_\theta(Y_t)\) and comparison with
\eqref{eq:initial-Jacobian}, started at \((r,Y_r)\), gives exactly
\eqref{eq:response-kernel}.  The bound follows from
\eqref{eq:uniform-flow-Jacobian}.
\end{proof}

\subsection{Gaussian models and Malliavin derivatives}

Let \(a:[0,T]\to[0,\Lambda]\) be deterministic and measurable.  On a
Wiener space, define
\begin{equation}\label{eq:deterministic-Gaussian-driver}
  X_t^a=\int_0^t\sqrt{a_r}\dd W_r,
\end{equation}
and let \(P^a\) be its law on \(\Omega\).  Then
\(P^a\in\mathfrak M_\Lambda\).  If
\(k\in H_0^1([0,T])\), translation of \(W\) by \(\eps k\) translates
\(X^a\) by the finite-variation path
\begin{equation}\label{eq:CM-driver-shift}
  h_t^{a,k}
  =\int_0^t\sqrt{a_r}\dot k_r\dd r.
\end{equation}

\begin{theorem}[Malliavin identification]\label{thm:malliavin-identification}
For every fixed \((\theta,s,t,y)\) and every \(1\le p<\infty\),
\[
  S_{s,t}^\theta(X^a;y)\in\mathbb D^{1,p}.
\]
Its Malliavin derivative with respect to the Brownian representation
\eqref{eq:deterministic-Gaussian-driver} is
\begin{equation}\label{eq:Malliavin-kernel}
  D_r^W S_{s,t}^\theta(X^a;y)
  =
  \one_{[s,t]}(r)\sqrt{a_r}\,
  \mathcal K_{s,t}^\theta(r;X^a,y)
\end{equation}
for almost every \(r\), almost surely.  In particular,
\begin{equation}\label{eq:Malliavin-uniform-bound}
  \norm{D^WS_{s,t}^\theta(X^a;y)}_{L^2([0,T])}
  \le
  \overline\sigma C_*
  \left(\int_s^t a_r\dd r\right)^{1/2}.
\end{equation}
\end{theorem}

\begin{proof}
Because the state derivatives of \(b_\theta\) and \(\sigma_\theta\) are
bounded, the classical Malliavin variational theorem for
\eqref{eq:modelwise-SDE} gives membership in \(\mathbb D^{1,p}\) for every
\(p>1\); see \cite[Theorem~2.2.1]{Nualart06}.  The case \(p=1\) follows
from \(\mathbb D^{1,p}\subset\mathbb D^{1,1}\) on the underlying probability
space.  For \(k\in H_0^1\), a Cameron--Martin shift of \(W\)
translates \(X^a\) by \(h^{a,k}\).  Moreover, the pathwise bound
\eqref{eq:response-kernel-bound} and the mean-value theorem in the shift
parameter give
\[
  \left|
    \frac{
      S_{s,t}^\theta(X^a+\eps h^{a,k};y)
      -S_{s,t}^\theta(X^a;y)
    }{\eps}
  \right|
  \le
  \overline\sigma C_*\Var(h^{a,k};[s,t])
  \le
  \overline\sigma C_*
  \left(\int_s^t a_r\dd r\right)^{1/2}
  \norm{\dot k}_{L^2}.
\]
The difference quotient therefore converges in every \(L^p\), by
\Cref{thm:bv-derivative} and dominated convergence.  The
Cameron--Martin shift characterization of the Sobolev derivative identifies
this limit with its Malliavin directional derivative and gives
\begin{align*}
  \langle D^WS_{s,t}^\theta,\dot k\rangle_{L^2}
  &=D_{h^{a,k}}S_{s,t}^\theta(X^a;y)\\
  &=\int_s^t
    \mathcal K_{s,t}^\theta(r;X^a,y)
    \sqrt{a_r}\dot k_r\dd r.
\end{align*}
Choose a countable dense subset \(\mathcal H_0\subset H_0^1([0,T])\).
Intersecting the probability-one sets above over \(k\in\mathcal H_0\), the
displayed identity holds simultaneously on one full set.  Both sides are
continuous linear functionals of \(\dot k\in L^2([0,T])\): for the right
side this follows from \eqref{eq:response-kernel-bound}.  Density extends
the identity to every Cameron--Martin direction and identifies the
\(L^2([0,T])\)-valued Malliavin derivative, proving
\eqref{eq:Malliavin-kernel}.  The uniform bound is
\eqref{eq:response-kernel-bound}.
\end{proof}

\begin{remark}[Scope of the Gaussian identification]
The pathwise derivative in \Cref{thm:bv-derivative} is defined without a
probability law.  Its interpretation as a Malliavin derivative requires a
specified Gaussian realization and Cameron--Martin space.  If the volatility
is deterministic, the Brownian Malliavin derivative is exactly the weighted
kernel in \eqref{eq:Malliavin-kernel}.  If the volatility is an adapted
functional of the Brownian path, a Cameron--Martin shift also differentiates
that control, and the corresponding variational equation contains additional
terms.
\end{remark}

\section{Conclusion}

Regularized dyadic square sums produce a Borel causal bracket and common
compact continuity cores for the nondominated class
\(\mathfrak M_\Lambda\).  The bracket drives a scalar \(C^1\) It\^o flow,
its initial Jacobian, and a coefficient-parameter tangent on the raw path
space.  Capacity-core transfer yields uniform \(L^p\) expansions,
Wasserstein stability of the first-jet laws, attainment for polynomial-growth
field payoffs, and a joint parameter--initial-state Danskin formula.
The constant-volatility example makes the corresponding active-model switch
explicit.
Finite-variation invariance gives a driver-response kernel whose Gaussian
realization is the Malliavin derivative after volatility weighting.  These
law-level conclusions therefore arise from one common causal analytic
construction.

{

}

\end{document}